\documentclass[]{interact}

\usepackage{epstopdf}
\usepackage[caption=false]{subfig}
\usepackage{array}
\newcolumntype{C}[1]{>{\centering\let\newline\\\arraybackslash\hspace{0pt}}m{#1}}
\usepackage{caption}
\usepackage{multirow}
\usepackage{bigstrut}
\usepackage{pdflscape}
\usepackage{setspace}
\usepackage[natbibapa,nodoi]{apacite}
\usepackage{hyperref}
\usepackage{amsmath}
\usepackage{mathtools}
\usepackage{lipsum}
\usepackage{color}
\usepackage{xcolor}
\usepackage{comment}
\usepackage{todonotes}
\usepackage[longnamesfirst,sort]{natbib}
\bibpunct[, ]{(}{)}{;}{a}{,}{,}
\renewcommand\bibfont{\fontsize{10}{12}\selectfont}

\theoremstyle{plain}

\theoremstyle{definition}

\theoremstyle{remark}

\begin{document}


\title{Optimal Minimal-Contact Routing of Randomly Arriving Agents through Connected Networks}

\author{
\name{Diptangshu Sen\textsuperscript{a}, Varun Ramamohan\textsuperscript{a}, and Prasanna Ramamoorthy\textsuperscript{b}\thanks{CONTACT Varun Ramamohan. Email: varunr@mech.iitd.ac.in} }
\affil{\textsuperscript{a}Department of Mechanical Engineering, Indian Institute of Technology Delhi, Hauz Khas, New Delhi 110016, India.}
\affil{\textsuperscript{b}Department of Management Studies, Indian Institute of Technology Delhi, Hauz Khas, New Delhi 110016, India.}
}

\maketitle

\begin{abstract}
Collision-free or contact-free routing through connected networks has been actively studied in the industrial automation and manufacturing context. Contact-free routing of personnel through connected networks (e.g., factories, retail warehouses) may also be required in the COVID-19 context. In this context, we present an optimization framework for identifying routes through a connected network that eliminate or minimize contacts between randomly arriving agents needing to visit a subset of nodes in the network in minimal time. We simulate the agent arrival and network traversal process, and introduce stochasticity in travel speeds, node dwell times, and compliance with assigned routes. We present two optimization formulations for generating optimal routes - no-contact and minimal-contact - on a real-time basis for each agent arriving to the network given the route information of other agents already in the network. We generate results for the time-average number of contacts and normalized time spent in the network. 
\end{abstract}


\section{INTRODUCTION}\label{intro}
Collision-free or contact-free routing of agents needing to navigate a connected network has been the subject of active research in recent years. This problem becomes particularly important in the manufacturing and logistics industry where a set of automated guided vehicles (AGVs) or robots are deployed in, for example, an assembly line manufacturing layout or a goods warehouse, to perform a set of tasks. For example, \cite{spensieri2015} consider the problem of planning the collision-free movement and scheduling of welding tasks performed in an automotive assembly line, and \cite{xin2020time} consider the problem of collision-free path planning and routing of multiple robotic effectors at an assembly line. In this work, we consider a related, but different problem: the real-time collision-free or contact-free route generation for an agent arriving randomly to a connected network, who is tasked with visiting a subset of nodes of the network in minimal time.

We define a contact as a non-zero overlap between the time intervals during which two or more agents are present at a network node. This implies that we consider contacts only at the nodes, and do not consider contacts between agents on the paths between nodes - we assume the paths have sufficient space to preclude catastrophic collisions or contacts between agents. Further, given the availability of free-ranging AGVs (i.e., AGVs that do not require preinstalled guide paths), collision avoidance along network edges is less likely to be a cause for concern \cite{duinkerken2006comparison,xin2020time}. The problem of generating such a contact-free route for an agent needing to visit a subset of nodes (hereafter referred to as the agent's \textit{node set}) in minimal time when a set of agents are already present in the network may arise in the context of AGVs deployed - one by one as the demand arises - to retrieve or deliver items from/to a set of stock points in a goods warehouse or machining stations in a factory layout. Similarly, this situation may also arise in the context of the COVID-19 pandemic. For example, it may be desirable to minimize contacts between human workers in a manufacturing/logistics context tasked with visiting and performing tasks at a set of destinations so as to minimize risk of infection transmission at said destinations. In this paper, we develop an optimization framework for assigning minimal contact routes to such agents and evaluate the framework - in terms of the number of contacts between agents and their time spent in the network - using a simulation wherein stochasticity is introduced in various aspects of the network traversal process of agents.

Our approach first involves constructing a simulation of the agent arrival and network traversal process. We then formulate the problem of assigning routes to each agent such that contacts are eliminated and/or minimized and shopping time is minimized as extensions of the Miller-Tucker-Zemlin (MTZ) integer programming formulation of the traveling salesman problem (TSP) \cite{miller1960}. The formulations are deterministic, and are parameterized using the expected values of the agent network traversal parameters, such as agent speed of movement from node to node, and the dwell time of an agent at each node. 

Given that we consider potential applications involving routing of human agents through a connected network, considering the impact of variation in the speed of agent movement and node dwell times (which may be variable for robotic agents as well) on the number of contacts becomes important. Further, compliance with the assigned no-/minimal-contact routes may also be of concern, and hence we model scenarios where agents can choose to deviate from their assigned routes in one of three ways. The need for a simulation to determine how these different forms of stochasticity in the network traversal process interact in terms of the number of contacts between agents thus becomes evident.

We now discuss the relevant literature and research contributions of our work.
We focus our survey of the literature on extensions of TSP and multi-agent routing problems that attempt to minimize contacts or collisions. There appears to be limited work in contact-free or collision-free routing from a stochastic simulation standpoint, likely because the majority of the work in this field has been carried out regarding collision-free routing for robotic systems, implying an assumption in these studies that uncertainty in network traversal aspects is likely to be minimal.    

A variant of the TSP similar to our current work is the $m$-TSP, wherein $m$ agents start and end their tour at the same node after visiting a set of nodes \cite{tolga2006}. Here, each node can be visited by only one agent. The objective of an $m$-TSP is to minimize the total tour length of all agents. Our problem differs from this problem in two respects: i) In $m$-TSP, all agents start out at the same time, while in our proposed formulations, the agents start at their simulated arrival times. ii) Second, in $m$-TSP each node is visited by only one agent, whereas in our case a particular node can be in the node set of multiple agents. More comprehensive discussions of the TSP and related variants can be found in the reviews by \cite{reinelt1991,applegate2006,cook2011,toth2014}.

As mentioned earlier, collision-free path planning of multiple robots or collision avoidance in AGV systems has been an active area of research. In many studies, the objective is to minimize the total cycle time of performing a set of tasks to be assigned to a set of robots/AGVs. Given a set of tasks to be performed by a set of agents, three problems are of interest: i) task allocation to multiple robots in the station \cite{herrero2010,korsah2013,lee2014}; (ii) routing problems to decide on the sequence in which the tasks are to be completed by each robot \cite{bullo2011,spensieri2015}; (iii) path planning to decide on the path taken by the robots while accomplishing their tasks \cite{guillaume2017}. Recent studies combine both routing and path planning \cite{xin2015control,miyamoto2016local,xin2020time}, and incorporate collision-free constraints in their analyses. Most studies in this area consider a static problem, as described above. A few studies consider dynamic problems, in terms of either the dynamic arrival of tasks \cite{duinkerken2006comparison,smolic2009time} or, to a limited extent, a dynamic environment \cite{duinkerken2006comparison,shi2018collision}. \cite{duinkerken2006comparison} - the only study we identified that explicitly used stochastic simulation in this context - explore routing strategies for free-ranging AGVs to avoid collisions with obstacles that arise on a stochastic basis; \cite{shi2018collision} also use simulation - for deterministically arriving agents - in a similar manner to evaluate their algorithm to find a feasible collision-free path. \cite{smolic2009time} consider finding the shortest path for an AGV to a single destination for dynamically arriving tasks, and avoid collisions by examining time windows in potential paths when other agents cross potential paths for the AGV under consideration. However, they do not simulate the dynamic arrival of tasks itself.

In relation to the above studies, our contributions are as follows. (i) We consider real-time generation of no-/minimal-contact routes for randomly arriving agents with independent sets of tasks (i.e., node sets). (ii) We propose generating no-/minimal-contact routes via an extension of the MTZ formulation of the TSP, yielding an exact solution. (iii) We explicitly consider the impact of uncertainty in agent speed, node dwell time, and compliance with the assigned route on the number of contacts via simulation.

We now describe the development of the agent arrival and network traversal simulation.

\section{AGENT ARRIVAL AND NETWORK TRAVERSAL SIMULATION}
\label{simopt}
In this section, we describe the development of the simulation of the agent arrival and network traversal process.

\subsection{Network Layout Description}
\label{simgen}
Without loss of generality, we model the connected network as a rectangular layout composed of identical unit square grids. The network we model has 70 grids, with 10 rows consisting of 7 grids each. Half of these 70 grids are nodes which the agents can visit, and the rest represent aisles through which agents can traverse the grid (move through the store) and access the nodes. Visualizations of the network along with example agent routes are provided in Figure~\ref{store}. We develop this type of layout to represent, for example, a warehouse with rows of stock points modelled as rows (or columns) of nodes and the grids between the nodes representing aisles through which a warehouse worker or an AGV can move from one stock point (node) to another.

The nodes in the network are categorized into different sections, keeping with the warehouse analogy where each section of a warehouse may contain a particular class of items. 
There is a fixed entry and exit node - in the simulation, the node [1,1] in Figure~\ref{store} represents both the entry and exit node. Based on an agent's node set, a feasible tour is one that starts at the entry node, visits all nodes in the agent's set at least once, and exits. 

\subsection{Simulation Details}

Even though simulating the stochastic arrival of agents is typical of discrete event simulation, our implementation of agent arrival and network traversal is more in line with a Monte Carlo simulation. 

We first generate a set of interarrival times for agents until the time of arrival of an agent exceeds the time horizon of the simulation. Once the arrivals of this set of agents are generated, we assign node sets to each agent (more details regarding this process is provided below). Routes to be assigned to each agent are generated based on their node set, their expected traversal speeds, expected node dwell times, and traversal pattern. The different types of traversal patterns we simulate for agents are described subsequently, including those based on the proposed TSP-NC and TSP-MC formulations. Stochasticity in traversal is introduced at three levels: (a) speed of movement of agents between nodes, (b) amount of time spent at each node (node dwell time), and (c) the compliance of an agent with the assigned route. Once a route is assigned to an agent, the set of actual speeds of movement between each pair of nodes on the route and the set of dwell times at each node are sampled from their respective distributions. Using the above set of information, the entry and exit times at each node are recorded for every agent. Contacts are determined by counting the number of overlaps between the entry and exit times of agents visiting each node. 

We now describe each aspect of the simulation in greater detail.

\textit{Arrival of Agents} : We assume that agent arrival follows a Poisson distribution, i.e., the inter-arrival times are exponentially distributed. For the rest of the analysis, we take the mean of the Poisson arrival distribution as $\lambda = 40$ agents/hour.

\textit{Node Set Generation for an Agent} : We consider a network traversal scenario where each agent has a predetermined set of nodes to visit. In a warehouse situation, this can represent a set of items that a worker or an AGV has to retrieve from various nodes in the warehouse. Each node is associated with a probability of featuring in an agent's node set. In a practical scenario (e.g., warehouse), this probability can be estimated from historical data regarding the frequency with which different nodes are visited. We set $p = 0.3$ for all nodes.  

\textit{Traversal Patterns}: In the real world, warehouse workers or AGV programmers may have a variety of ways in which they traverse their node set. In order to capture this diversity of traversal patterns, we assume (and simulate) two principal traversal patterns against which we benchmark our optimal routes.

The first type is a \textit{greedy} traversal pattern where an agent chooses which section to visit and subsequently how to traverse all required nodes in that section, greedily. Upon entry, the said agent quickly identifies which section is the closest of all sections he needs to visit and proceeds there. Within the section, he also follows a greedy approach, moving to the next nearest node until he has visited all required nodes in that section. The next section is again chosen in a greedy manner, based on which of the remaining sections are closest from his current position.

The second type is a \textit{preferential} traversal pattern, based on agents who have a predetermined order in which they wish to visit different sections in the network. However, we assume that inside a section, they visit nodes in a greedy manner, always moving to the next closest node. In a warehouse setting, this can represent a worker choosing to visit a section with heavier items that can be loaded on the bottom of a forklift pallet and then visiting a section with lighter items that can be loaded on top of the heavier items on the pallet.

\textit{Node Dwell Time}: We assume that the time spent by an agent at a node is exponentially distributed with a mean of 2 minutes. We choose the exponential distribution as a worst-case scenario where the amount of time an agent has spent at a node does not provide any information regarding the completion of the dwell time of the agent at the node.

\textit{Agent Speed}: We assume that the speed of agent movement $v$ between nodes is a random variable given by $v = 10(1+x)$, where $x$ is a beta random variable with $\alpha =0.5$ and $\beta = 1.5$. $v$ is thus distributed in the range 10 meters/minute to 20 meters/minute, with a mean speed of 12.5 meters/minute. 

In addition to the above aspects of the simulation, we also introduce stochasticity in terms of the compliance of agents with assigned routes; however, we discuss this in detail after we present the TSP-NC and TSP-NC optimization formulations given that the notion of compliance becomes of interest when routes are assigned or recommended to an agent as opposed to when an agent chooses their own routes.

In Figure~\ref{store} below, we depict the paths of both agent traversal patterns. The grey grids represent nodes and the white grids represent paths (e.g., aisles in a warehouse) for movement for the agents. There are a total of nine sections available in the network, $a$ through $i$. Each grid is of size 5m $\times$ 5m. We use the same specifications of the network throughout our analysis.

In each case in Figure~\ref{store}, we assume the agent's node set consists of six nodes, corresponding to nodes at positions $(1,5), (3,2), (4,4), (6,2), (6,4)$ and $(8,4)$ in the network. For the preferential node set traversal pattern (the figure on the right), the order in which the agent prefers to visit the sections associated with their list is given by $B -> F -> G -> C -> A$.

 \begin{figure}[!ht]
    \begin{minipage}[b]{0.49\linewidth}
    \centering
    \includegraphics[width = 0.65\textwidth]{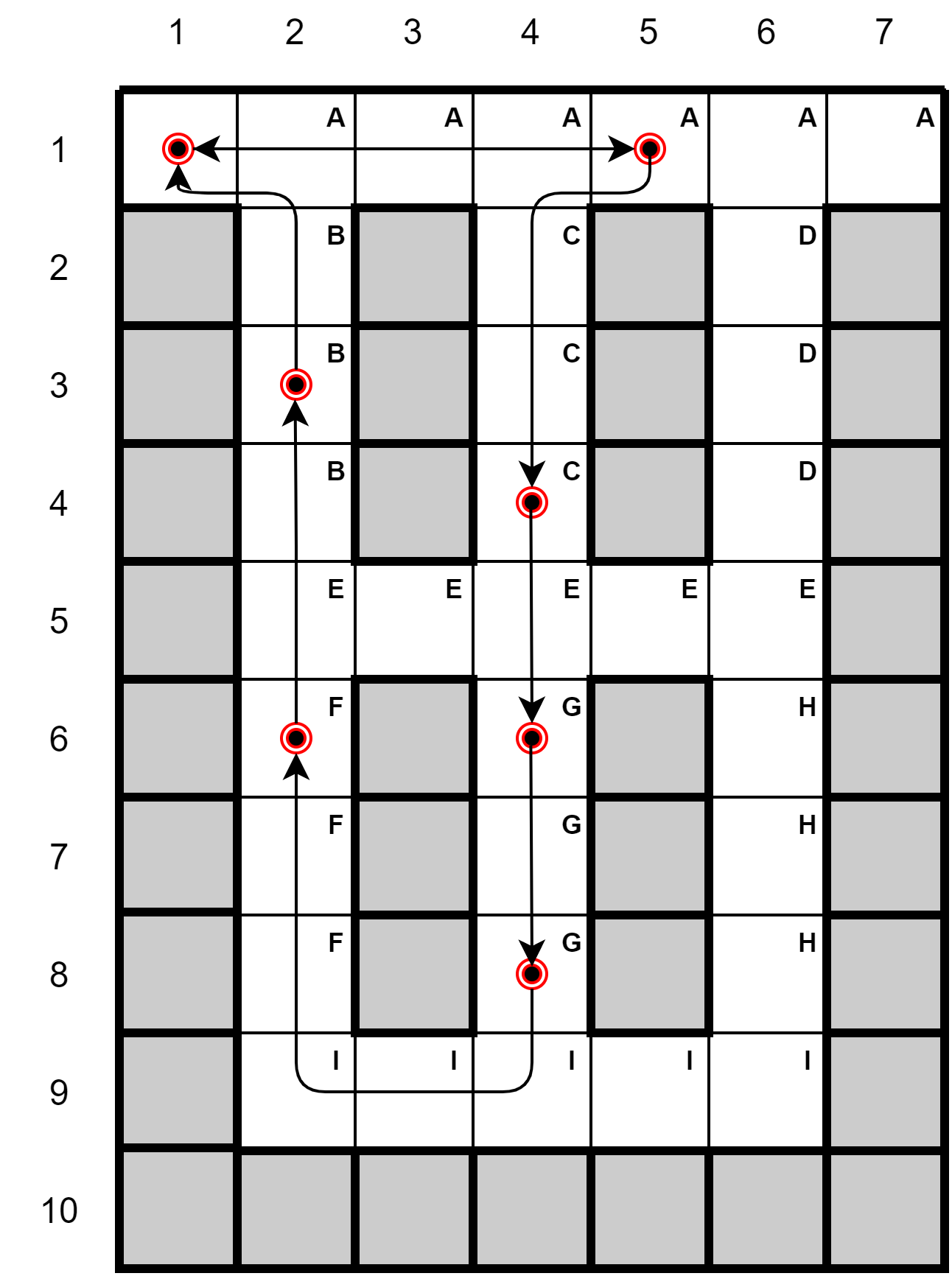}
    \end{minipage}
    \begin{minipage}[b]{0.49\linewidth}
    \centering
    \includegraphics[width = 0.65\textwidth]{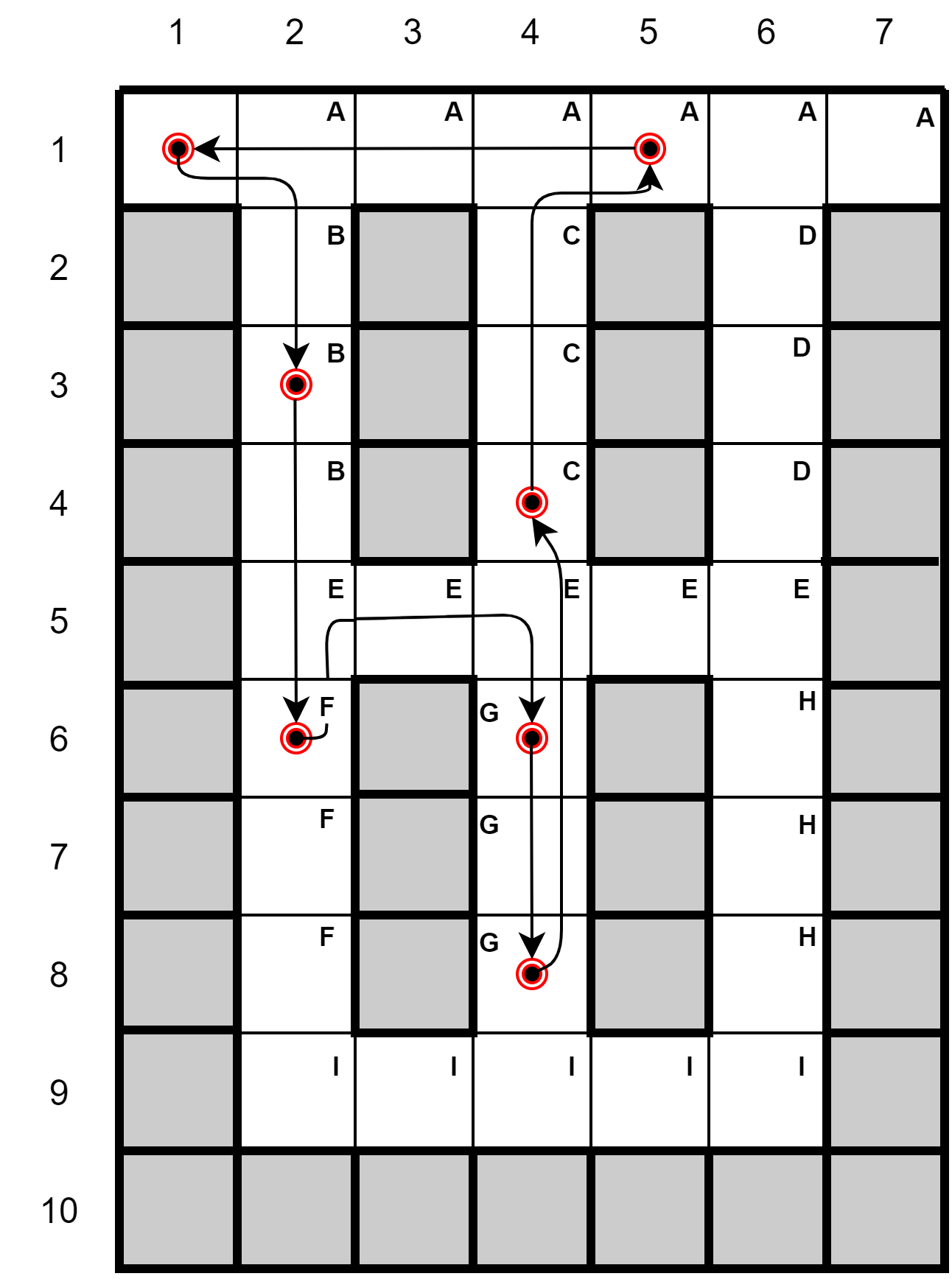}
    \end{minipage}\hfill
\caption{Example routes for greedy (left) and preferential (right) agents.}
\label{store}
\end{figure}

We present the simulation results for the greedy and preferential node set traversal patterns along with the results from the extended TSP formulations.

\section{NO-CONTACT/MINIMAL-CONTACT ROUTE GENERATION}
\label{tsp}
For a given agent arrival rate, the greedy and preferential traversal patterns lead to a large number of contacts among agents. Here a contact is defined as the event where at least two agents are present at the same node at the same time. Further, these traversal patterns are also non-optimal in terms of minimizing the time spent in the network. Hence we propose to generate a customized route for each agent, based on his node set, that eliminates or minimizes contacts with other agents and also minimizes the time spent in the network. Our optimization framework builds on the MTZ formulation of the TSP \cite{reinelt1991}. We present the no-contact formulation (referred to hereafter as the TSP-NC formulation hereafter) first.

\subsection{No-contact TSP Formulation}
\label{tspnc}

This formulation is a mixed-integer program that provides an exact solution. The formulation is based on recognizing that in order to avoid contacts, an agent should avoid reaching a node when it is already occupied by another agent. To implement this, we introduce a set of continuous time variables that capture the time points when an agent reaches the different nodes on their tour. These time variables are then constrained to avoid overlap with the time windows when the nodes are occupied/blocked by other agents. The blocked time windows for each node are provided as input data to the formulation based on the routes generated for agents already present in the network. The formulation is provided below.

\textit{Parameters}
\begin{itemize}
    \item $c_{ij}$: Length of shortest path joining nodes $i$ and $j$. Here $i, j \in \{1,2,...,n\}$, where $n$ is the number of nodes.
    \item $E[v]$: Expected node-to-node speed of movement for an agent.
    \item $E[b]$: Expected node dwell time.
    \item $d_i = \{d_{i1}, d_{i2},...,d_{iK_i}\}$: A set of $K_i$ time points indicating the time points $d_{ik} \in d_i$ at which node $i$ will start to be occupied for the next $E[b]$ minutes by $K_i$ other agents in the network, $i \in N$. Note that $K_i$ will be dependent on each node $i$, and hence the subscript.
    \item $t_0$: The time instant when a new agent arrives.
    \item $M$: A large positive number.
\end{itemize}

\textit{Decision Variables}
\begin{itemize}
    \item $x_{ij}$ $\in$ \{0, 1\}: Boolean variable with value 1 if the tour includes travel from node $i$ to node $j$; 0 otherwise.
    \item $u_i$ $\in$ $\mathbb{Z}^{+} \; \forall \; i \in N$: Dummy variable indicating the index at which node $i$  features in the ordered tour.
    \item $t_i$ $\in$ $\mathbb{R}^{+}\cup\{0\} \; \forall \; i \in N$: Time at which node $i$ is reached.
    \item $y_{ik}$ $\in \{0, 1\} \; \forall \; i \in N, \;k \in K_i$: Dummy Boolean variable for either-or constraints.
\end{itemize}
\begin{equation*}
    \begin{aligned}
        &\text{min} \hspace{3pt} \sum_{i=1}^{n}\sum_{j=1}^{n} c_{ij}x_{ij} \hspace{20pt}subject\hspace{3pt}to:\\
        &\text{MTZ formulation constraints:}\\
        &\;\; \;\; x_{ii} = 0, \hspace{6pt}\forall\hspace{3pt} 1 \leq i \leq n\\
        &\;\; \;\; \sum_{i=1}^{n}x_{ij} = 1,\hspace{6pt}\forall\hspace{3pt}1 \leq j \leq n\\ 
        &\;\; \;\; \sum_{j=1}^{n}x_{ij} = 1,\hspace{6pt}\forall\hspace{3pt}1 \leq i \leq n\\
        & \;\; \;\; u_1 = 1\\ 
        &\;\; \;\; 2 \leq u_i \leq n, \;\forall \hspace{6pt}2\leq i\leq n\\
        &\;\; \;\; u_i - u_j +nx_{ij} \leq n-1, \;\forall\hspace{6pt} 1\leq i \leq n,\; 2 \leq j \leq n,\; i\neq j\\
        &\text{Time constraints:}\\
       & \;\; \;\; t_1 = t_0\\
       & \;\; \;\; t_j - t_i \leq M(1-x_{ij})+E[b]+\frac{c_{ij}}{E[v]},\hspace{6pt}\forall \hspace{3pt} 1\leq j\neq i \leq n\\
       & \;\; \;\; t_j - t_i \geq -M(1-x_{ij})+E[b]+\frac{c_{ij}}{E[v]},\hspace{6pt} \forall \hspace{3pt} 1\leq j\neq i \leq n\\
       \end{aligned}
\end{equation*}
\begin{equation*}
    \begin{aligned}
       &\text{TSP-NC constraints:}\\
       & \;\; \;\; t_i + E[b] \leq d_{ik} + My_{ik}, \hspace{6pt}\forall\hspace{3pt} 2 \leq i \leq n, \; \; 1\leq k \leq K_i;\\    
       & \;\; \;\; t_i + M(1-y_{ik}) \geq d_{ik} + E[b], \hspace{6pt}\forall\hspace{3pt}2 \leq i \leq n, \; \; 1\leq k \leq K_i\\
\end{aligned}
\end{equation*}
The first (bookkeeping) constraint prevents an edge looping into the same node. The second and third constraints allow only one edge to be directed towards and away from a given node, forcing the fact that each node is visited only once. $u_1 = 1$ and $t_1 = t_0$ initialize that the tour starts at the entry node at time $t_0$. The purpose of the constraints involving $u_i$ and $u_j$ is to eliminate subtours. The new $t_i$, $t_j$ constraints are in place to force $t_j = t_i + E[b] + \frac{c_{ij}}{E[v]}$ if an edge exists between nodes $i$ and $j$. If an edge does not exist between nodes $i$ and $j$, then the constraints ensure enough slack between $t_i$ and $t_j$ in either direction to allow for sufficiently large differences between $|t_i - t_j|$. The last two no-contact constraints ensure elimination of overlaps between the intervals when the current customer occupies the $i^{th}$ node (the interval $[t_i, t_i + b]$) with the interval $[d_{ik}, d_{ik}+b]$ during which node $i$ is expected to be occupied by another agent in the network. 

\begin{figure}[htbp]
     \centering
     \includegraphics[width = 0.85\textwidth]{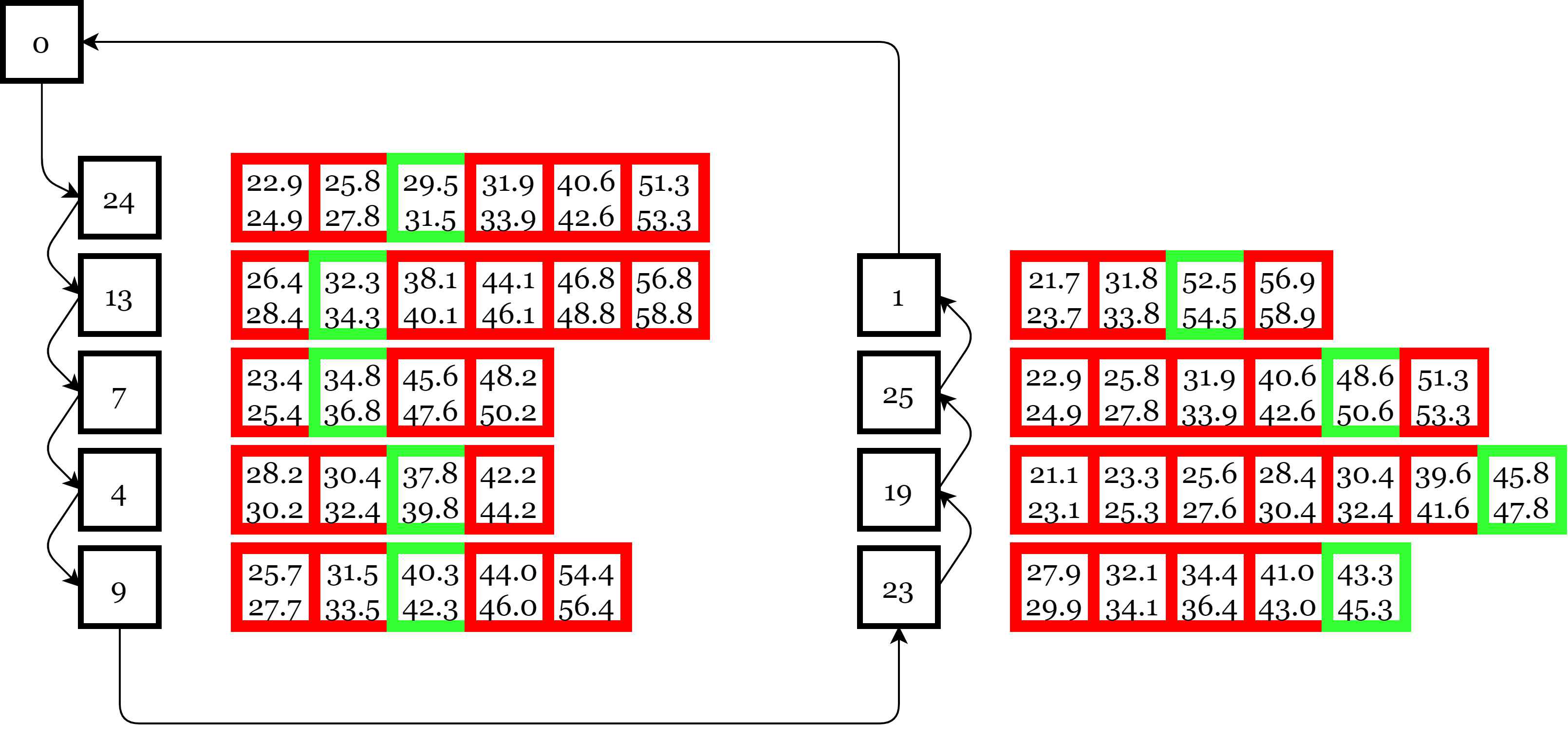}
     \caption{Depiction of a sample optimal no-contact route. The flowchart indicates the recommended order of traversal of nodes. For a given node, the red blocks denote the time windows during which the node is occupied while the green block denotes the assigned time window for visiting said node (as per TSP-NC formulation). Observe that there is no overlap with a blocked window, indicating zero contacts.}
     \label{tspncrt}
 \end{figure}

Thus the formulation accounts for the fact that there can be multiple time windows during which a node is blocked. The no-contact path generation exercise is carried out when a new agent arrives at the network with a predetermined set of nodes they wish to visit. Their assigned path avoids contacts with all agents already present in the network. A sample optimal route generated for an agent under this formulation is depicted in Figure~\ref{tspncrt}.

Finding an optimal or even a feasible no-contact solution within reasonable runtimes can be challenging under certain conditions. In some cases, for a new agent who has just arrived, all the nodes in their node set may be blocked such that an optimal no-contact route cannot be generated without a change in the state of the network (e.g., an agent leaving the network after finishing their tour). In the framework of the above formulation, this implies infeasibility. To address this, in the next subsection, we present a relaxation of the TSP-NC formulation that generates a route which is time-efficient and also minimizes contacts to the extent possible. 

\subsection{Minimal-contact Formulation}
\label{tspmc}
We refer to this relaxation of the TSP-NC formulation as the TSP-MC formulation (the `MC' stands for minimal contact). This formulation will yield a route with zero contacts if such a route is feasible and yields a lower objective function value when compared to a route with one or more contacts. However, if such a route is not feasible, it will yield an optimal route with as few contacts as possible.

For this formulation, in addition to the parameters and decision variables introduced for the TSP-NC formulation, we introduce the penalty parameter $P_e$, to be applied each time a contact occurs between agents. We also introduce two new decision variables as given below. The variables $y_{ik}$ from the TSP-NC formulation are not used in the TSP-MC formulation below. In the formulation below, due to space limitations, we only present the objective function and the minimal-contact constraints. The MTZ TSP and time constraints remain the same as in the TSP-NC formulation, and the non-negativity constraints follow from the definition of the decision variables.

\begin{itemize}
    \item $\delta_{ik}$ $\in$ $\{0, 1\}, \; \forall \; i \in N, k \in K_i$ :  Boolean variable that records contacts between a pair of agents (i.e., $\delta_{ik} = 1$ if a contact occurs between two agents at node $i$, 0 otherwise).
    \item $z_{ik}$ : continuous non-negative variable used to determine whether contacts between the current agent and another agent occur at node $i$ based on the absolute value of the difference between $t_i$ and $d_{ik}$. The value of $z_{ik}$ is determined by the value of $t_i$.
\end{itemize}
\begin{equation*}
    \begin{aligned}
        &\text {min} \hspace{3pt} \sum_{i=1}^{n}\sum_{j=1}^{n} c_{ij}x_{ij} + P_e \sum_{i=1}^{n}\sum_{k = 1}^{K_i}\delta_{ik} \hspace{20pt} subject\hspace{3pt}to:\\
        &\text{TSP-MC constraints}:\\
        & \; \; z_{ik} = |(t_i-d_{ik})/b|, \hspace{6pt}\forall\hspace{3pt}1 \leq i \leq n, 1 \leq k \leq K_i\\
        & \; \; M\delta_{ik} + z_{ik} \geq 1, \hspace{6pt}\forall\hspace{3pt}1 \leq i \leq n, 1 \leq k \leq K_i\\
        & \; \; -M(1-\delta_{ik}) + z_{ik} \leq 1, \hspace{6pt}\forall\hspace{3pt}1 \leq i \leq n, 1 \leq k \leq K_i\\
    \end{aligned}
\end{equation*}

 The first constraint, via the variable $z_{ik}$, records the extent of overlap between the time spent by the current agent at the $i^{th}$ node and the time spent by any other agent at the same node. If $z_{ik} > 1$, then there is no overlap, and hence no contact occurs; otherwise there is a contact. The last two TSP-MC constraints ensure that if a contact occurs, then the $\delta_{ik}$ variable is set to 1, thus recording the contact.

Note that $P_e$ can serve as a tuning parameter for the network administrator. A higher value of the penalty lowers the number of contacts at the expense of time-efficiency.


\subsection{Noncompliance with Assigned TSP-NC/TSP-MC Routes}
\label{noncomp}
The TSP-NC and the TSP-MC routes are assigned to agents based on the expected values of the node-to-node speed and node dwell time random variables. Therefore, the number of contacts recorded during the simulation due to uncertainty in these parameters will be higher than those promised by the formulations. Another source of increased contacts when deploying these formulations in practice will involve non-compliance on part of the agents with assigned TSP-NC/TSP-MC routes. We describe our approach towards modelling noncompliance with the assigned routes now.

The issue of noncompliance may be particularly relevant where the agents are human workers navigating a network such as a warehouse; in the case of AGVs, this is likely to be less of an issue, unless they are under real-time control by a human programmer. We consider three types of noncompliance. 

\textit{Type A noncompliance}. This type of noncompliance can occur in the following situation. Consider an agent currently at node $i$ in their assigned route, with node $j$ being next in their assigned route. We refer to this node as the \textit{next assigned} node. Let another node $k$ (distinct from node $j$), also in the node set of the agent, be located at a smaller distance from node $i$ (we refer to this node as the \textit{greedy} node). In this situation, an agent may deviate, with probability $p_{ijk}$, from their assigned route by moving to the greedy node $k$ instead of the next assigned node $j$. We assume that $p_{ijk}$ is a function of the distances between node $i$ and nodes $j$ and $k$ ($d_{ij}$ and $d_{ik}$, respectively). We model $p_{ijk}$ as follows:
    \begin{equation*}
   p_{ijk} =
\left\{
	\begin{array}{ll}
		0  & \mbox{if } d_{ij} \leq d_{ik} \\
		1 - \frac{d_{ik}}{d_{ij}} & \mbox{otherwise } 
	\end{array}
\right.
    \end{equation*}
Thus the probability that an agent is noncompliant with their assigned route when at a given node, depends on the extent of the difference in the distances between their current node and the greedy and next assigned nodes, respectively. Finally, we assume that once an agent engages in this type of noncompliance, they abandon the assigned route completely, and visit the remainder of their nodes in a greedy manner.

\textit{Type B noncompliance.} This type of noncompliance is similar to type A noncompliance; however, once the agent deviates to the greedy node, they return to the TSP-NC/TSP-MC route by moving back to the originally next assigned node from the greedy node. Subsequently, they follow the assigned route, skipping over the greedy node (that was visited out of turn) when it appears in the route. Note that once the agent returns to the assigned route from the greedy node, they can still engage in noncompliance with regard to each next assigned node with probability $p_{ijk}$.

\textit{Type C noncompliance}. This type of noncompliance is different from types A and B. At each node, an agent may choose to visit one of the neighboring nodes \textit{not} in their node set with a certain probability (equal for each of a node's neighbors). For the network under consideration, because of its rectangular grid-like structure, we define the neighbors of a node as the set of nodes that share an edge (as opposed to a vertex) with said node. At an assigned node, the agent may visit any one of its neighbors \textit{not} in its node set with a total probability equal to 0.2. Once the agent visits a neighbor, the agent does not return to the assigned node; instead, they directly move to the next node on the assigned route. We incorporate this type of noncompliance to represent situations wherein a worker (human or AGV) in a warehouse may receive impromptu instructions to retrieve an item from a neighboring stock point that was not originally in their node set. 

\begin{figure}[htb]
	\centering
	\begin{minipage}[b]{0.48\linewidth}
		\centering
		\includegraphics[width=0.90\textwidth]{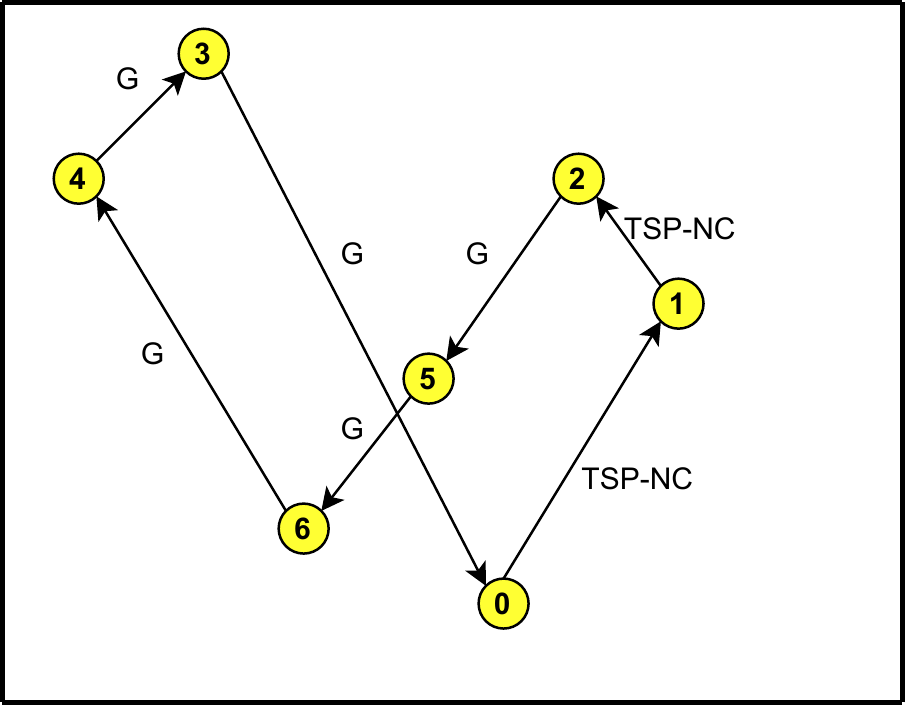}
		\caption{Type A noncompliance.}
		\label{typea}
	\end{minipage}
	\begin{minipage}[b]{0.48\linewidth}
		\centering
		\includegraphics[width=0.90\textwidth]{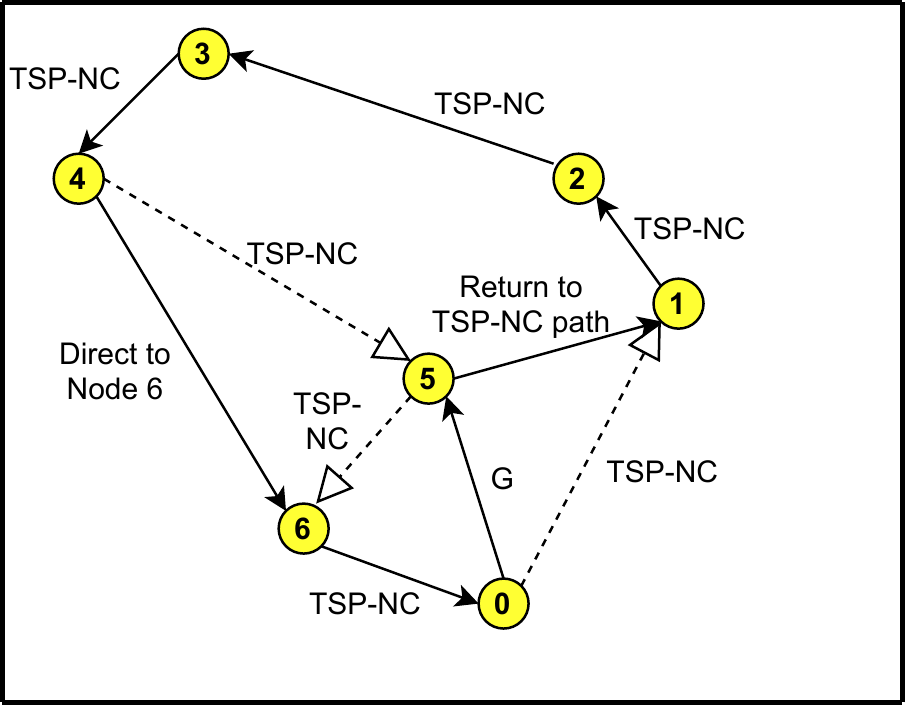}
		\caption{Type B noncompliance.}
		\label{typeb}
	\end{minipage}%
	\caption{Type A and B noncompliance examples. Solid lines in the figure denote the actual route taken by the agent; dotted lines denote the assigned TSP-NC path. `G' denotes the greedy path taken at that node.}
	\label{fignoncomp}
\end{figure}

Noncompliance of types A and B are depicted in Figure~\ref{fignoncomp} below (we do not depict Type C noncompliance due to space limitations). We now present the numerical results from implementing these formulations within the simulation.

\section{SIMULATION OUTCOMES}
\label{numres}
In Table~\ref{nononcomp}, we present summary statistics for relevant outcomes - the number of contacts and mean normalized time spent in the network for an agent - for all node traversal patterns when noncompliance is not involved. The mean normalized network time (MNNT) for an agent is calculated as the ratio of the total time spent in the network for an agent to the cardinality of its node set. Under this scenario (no noncompliance), we explore two possible sub-scenarios. The first involves a situation where only agent arrivals to the network are random, and the speed of movement and node dwell times are deterministic. This can represent a scenario wherein the agents in question are AGVs deployed to retrieve a set of items from a warehouse, as AGVs are likely to have nearly constant movement speeds and item retrieval times at stock points. The second sub-scenario involves incorporating uncertainty in agent speed of movement as well as node dwell times, which can model a scenario wherein workers in a warehouse have strong disincentives for deviating from the assigned route, but their speeds of movement and dwell times at a stock point are variable. We do not consider a scenario wherein uncertainty is incorporated in the agent speed of movement alone or in the node dwell time alone because it is unlikely that uncertainty in only one of the parameters will be present. Further, even if uncertainty is present in only one of these parameters (but not both), this scenario can easily be incorporated within our framework.

The results in Table~\ref{nononcomp} are generated based on three hours of simulation time (i.e., three hours of network operation), wherein the first one hour is considered as a warm-up period, and outcomes are collected for agents arriving from the beginning of the second hour onwards (i.e., the number of contacts reported are those observed over these two hours). The summary statistics are obtained from 30 replications for each node traversal pattern. The simulation was programmed in Python, and the Gurobi optimization suite \cite{gurobi} was used to solve the TSP-NC/-MC formulations. The majority of the numerical experiments were performed on an Intel Xeon workstation with 4 cores and 32 gigabytes of memory. The MNNT reported in Table~\ref{nononcomp} is associated with the case when only arrivals are random; we do not report the MNNT for the scenario with random agent arrivals, node dwell time, and agent speed as it remains approximately the same as that observed in the scenario with only random agent arrivals. 

\begin{table}[htbp]
  \centering
  \caption{Simulation outcomes without considering noncompliance. The mean and standard deviation of the number of contacts are reported; SE = standard error of the mean.}
    \begin{tabular}{|C{2.5cm}|C{3.5cm}|C{3cm}|C{4.5cm}|}
    \hline
    Node traversal pattern & Contacts (random arrivals) & MNNT (SE) & Contacts (random arrivals, dwell time, speed) \bigstrut\\
    \hline
    Greedy & 740.10 (192.72) & 3.51 (0.04) & 660.07 (126.37) \bigstrut\\
    \hline
    Preferential & 697.10 (133.85) & 4.01 (0.04) & 706.8 (144.84) \bigstrut\\
    \hline
    TSP   & 992.1 (280.60) & 3.15 (0.04) & Not relevant \bigstrut\\
    \hline
    TSP-NC & 0 & 3.25 (0.08) & 293.3 (60.79) \bigstrut\\
    \hline
    TSP-MC & 0 & 3.36 (0.04) & 293.7 (59.61) \bigstrut\\
    \hline
    \end{tabular}%
  \label{nononcomp}%
\end{table}%

It is evident from Table~\ref{nononcomp} that the TSP-NC and TSP-MC routes yield a significant reduction in the number of contacts when compared to the no-/minimal-contact node traversal patterns. This is particularly evident when only stochastic arrivals are considered, and hence this indicates that this approach could be particularly relevant when variation in node dwell times and agent speeds is minimal. Even with stochastic dwell times and agent speeds, we see that the reduction in the number of contacts is substantial. Interestingly, with the greedy and preferential node traversal patterns, we see that the number of contacts do not increase with stochastic node dwell times and agent speeds in comparison to the situation with only stochastic agent arrivals; this implies an effect akin to `destructive interference' between the stochasticity in different parameters.

With regard to MNNT, it can be seen that the TSP-NC and TSP-MC routes yield reductions between 4.27\% - 18.95\% when compared to the greedy and preferential node traversal patterns. Note that the increase in MNNT of the TSP-NC pattern when compared to the MTZ TSP without no-/minimal-contact constraints is not statistically significant at a 5\% level of significance. In comparison to the TSP-MC pattern, the increase of 6.67\% is statistically significant at a 5\% level of significance.

The mean computational runtimes of the TSP-NC and TSP-MC models for generating no-/minimal-contact routes are 55.82 (standard error [SE] = 22.91) seconds, and 3.88 (SE = 3.26) seconds. The runtimes for the TSP-NC formulation vary widely depending upon: (a) the cardinality of the agent's node set, (b) the average cardinalities of the node sets of the agents that are taken into account for the generation of the TSP-NC route for the agent under consideration node sets, and (c) the extent to which the nodes in the node set of the agent under consideration are occupied by other agents during the time likely to be spent in the network by the agent under consideration. These conditions can even lead to infeasibilities when the TSP-NC formulation is used to generate routes for arriving agents. However, this occurs infrequently, and we observe an average number of 4.5 (SD = 2.88) infeasibilities (out of approximately 80 agents) for the TSP-NC formulation for the two cases considered in Table~\ref{nononcomp}. Note that the number of infeasibilities is not affected by the extent of stochasticity (e.g., in agent arrivals only, or in agent speed and node dwell times as well) in the simulation as the TSP-NC/-MC formulations are solved using expected values of these parameters.

We also note that in the stochastic arrivals only scenario the TSP-MC formulation also yields zero contacts but with substantially lesser runtimes and a marginally higher MNNT. A similar trend is also observed with stochastic node dwell times and agent speeds, with the average number of contacts nearly equal to that under the TSP-NC node traversal pattern. This indicates that the TSP-MC formulation may be particularly suited for this purpose; however, the extent of the increase in the number of contacts when the agent arrival rate increases must be ascertained prior to choosing one formulation over the other if only a single formulation can be deployed in a practical setting.

We present the results for a scenario where all types of noncompliance are incorporated along with stochasticity in agent arrivals, node dwell times, and in agent speed. This scenario may also apply to the warehouse setting with human workers or in a setting with real-time programming and/or control of AGVs by human operators. We present the results in Table~\ref{comptab} for a scenario where each agent has an equal probability of engaging in type A, B or C noncompliance.

\begin{table}[htbp]
  \centering
  \caption{Simulation outcomes with agent noncompliance with assigned routes. SD = standard deviation; SE = standard error.}
    \begin{tabular}{|C{4cm}|c|C{3cm}|}
    \hline
    Node traversal pattern & Contacts: mean (SD) & MNNT (SE) \bigstrut\\
    \hline
    TSP-NC & 390.7 (61.1) & 3.58 (0.09) \bigstrut\\
    \hline
    TSP-MC & 344.9 (84.51) & 3.58 (0.09) \bigstrut\\
    \hline
    \end{tabular}%
  \label{comptab}%
\end{table}%

From Table~\ref{comptab}, it is clear that introducing noncompliance increases both the number of contacts and the MNNTs; however, the average number of contacts with the TSP-NC/-MC routes are still substantially lesser when compared to greedy and preferential node traversal patterns. Reductions in the average number of contacts range from 41-50\%. This result can be considered to a `worst-case' scenario, given that we assume that every agent engages in noncompliance; further reductions in the average number of contacts are possible if this assumption is relaxed. 

We also see from the above results that generating TSP-NC/-MC routes for agents also yields a substantial reduction in the variances of the number of contacts, even with noncompliance and stochasticity introduced in all possible parameters associated with agent node traversal. 

The above results are generated for a case with the probability of an agent engaging in a type of noncompliance being equal for all types; however, all types of noncompliance do not appear to yield the same outcome in terms of the average number of contacts. Simulation experiments performed with stochastic arrivals and agents engaging in only one type of noncompliance yielded the following number of contacts: (a) 209.2 (SD = 36.17) contacts with type A; (b) 314.7 (SD = 78.65) contacts with type B, and (c) 373.3 (SD = 125.49) contacts with type C noncompliance. Note that the variance in the number of contacts also changes with the type of noncompliance. Hence instituting stronger disincentives to engage in one type of noncompliance over another may be worth considering for the network administrator.

We now conclude the paper with a summary of the work, its potential impact and a discussion of study limitations.

\section{CONCLUSIONS \& DISCUSSION}
\label{discuss}
In this study, we present optimization formulations for assigning contact-free routes to agents arriving randomly to a connected network. We evaluate - via a Monte Carlo simulation - the efficiency of the assigned routes in terms of the number of collisions/contacts between agents traversing their node sets in the network under stochasticity in various parameters associated with traversal. Under all scenarios, we demonstrate that our formulations yield significant reductions in contacts between agents with only marginal increases in shopping time when compared to multiple commonly used but sub-optimal (in terms of minimizing contacts) network traversal patterns.

Among the two formulations that we propose, it is evident that the TSP-MC formulation is computationally more efficient than the TSP-NC formulation while yielding the same or marginally higher number of contacts. However, we have presented results only for a single agent arrival rate in this paper. Preliminary experiments indicate an increase in the number of contacts with an increase in the arrival rate with the TSP-MC formulation (while remaining computationally efficient); however, the increase is marginal when compared to the overall jump in contacts that would result from higher arrival rate in the absence of an organized route-planning framework like ours. 

A limitation of this work is that the construction of the network, specification of traversal patterns, types of noncompliance, and the parameterization of the simulation and optimization framework is not based on data from the real-world; hence an immediate avenue of future research involves developing a version of the simulation and optimization framework using data from a real-world scenario. Further, the no-/minimal-contact optimization frameworks are based on extensions to the MTZ formulation of the TSP, and hence they may be inefficient when deployed for networks of substantially larger size than what is considered in this paper. Hence another direction of future research involves development of the TSP-NC/-MC equivalents with more efficient versions of the TSP, such as the TSP with time windows.

\section*{ACKNOWLEDGEMENTS}
The authors gratefully acknowledge S.G. Deshmukh and Subhangshu Sen for their useful suggestions during the course of the work.

\footnotesize

\bibliographystyle{apacite}
\bibliography{demobib}

\begin{thebibliography}{}

\bibitem [\protect \citeauthoryear {%
Applegate%
, Bixby%
, Chvatal%
\BCBL {}\ \BBA {} Cook%
}{%
Applegate%
\ \protect \BOthers {.}}{%
{\protect \APACyear {2006}}%
}]{%
applegate2006}
\APACinsertmetastar {%
applegate2006}%
\begin{APACrefauthors}%
Applegate, D\BPBI L.%
, Bixby, R\BPBI E.%
, Chvatal, V.%
\BCBL {}\ \BBA {} Cook, W\BPBI J.%
\end{APACrefauthors}%
\unskip\
\newblock
\APACrefYear{2006}.
\newblock
\APACrefbtitle {The Traveling Salesman Problem: A Computational Study} {The
  traveling salesman problem: A computational study}.
\newblock
\APACaddressPublisher{}{Princeton University Press}.
\PrintBackRefs{\CurrentBib}

\bibitem [\protect \citeauthoryear {%
Bektas%
}{%
Bektas%
}{%
{\protect \APACyear {2006}}%
}]{%
tolga2006}
\APACinsertmetastar {%
tolga2006}%
\begin{APACrefauthors}%
Bektas, T.%
\end{APACrefauthors}%
\unskip\
\newblock
\APACrefYearMonthDay{2006}{}{}.
\newblock
{\BBOQ}\APACrefatitle {The Multiple Traveling Salesman Problem: An Overview of
  Formulations and Solution Procedures} {The multiple traveling salesman
  problem: An overview of formulations and solution procedures}.{\BBCQ}
\newblock
\APACjournalVolNumPages{Omega}{34}{3}{209--219}.
\PrintBackRefs{\CurrentBib}

\bibitem [\protect \citeauthoryear {%
Bullo%
, Frazzoli%
, Pavone%
, Savla%
\BCBL {}\ \BBA {} Smith%
}{%
Bullo%
\ \protect \BOthers {.}}{%
{\protect \APACyear {2011}}%
}]{%
bullo2011}
\APACinsertmetastar {%
bullo2011}%
\begin{APACrefauthors}%
Bullo, F.%
, Frazzoli, E.%
, Pavone, M.%
, Savla, K.%
\BCBL {}\ \BBA {} Smith, S\BPBI L.%
\end{APACrefauthors}%
\unskip\
\newblock
\APACrefYearMonthDay{2011}{}{}.
\newblock
{\BBOQ}\APACrefatitle {Dynamic Vehicle Routing for Robotic Systems} {Dynamic
  vehicle routing for robotic systems}.{\BBCQ}
\newblock
\APACjournalVolNumPages{Proceedings of the IEEE}{99}{9}{1482--1504}.
\PrintBackRefs{\CurrentBib}

\bibitem [\protect \citeauthoryear {%
Cook%
}{%
Cook%
}{%
{\protect \APACyear {2011}}%
}]{%
cook2011}
\APACinsertmetastar {%
cook2011}%
\begin{APACrefauthors}%
Cook, W\BPBI J.%
\end{APACrefauthors}%
\unskip\
\newblock
\APACrefYear{2011}.
\newblock
\APACrefbtitle {In Pursuit of the Traveling Salesman: Mathematics at the Limits
  of Computation} {In pursuit of the traveling salesman: Mathematics at the
  limits of computation}.
\newblock
\APACaddressPublisher{}{Princeton University Press}.
\PrintBackRefs{\CurrentBib}

\bibitem [\protect \citeauthoryear {%
Duinkerken%
, Ottjes%
\BCBL {}\ \BBA {} Lodewijks%
}{%
Duinkerken%
\ \protect \BOthers {.}}{%
{\protect \APACyear {2006}}%
}]{%
duinkerken2006comparison}
\APACinsertmetastar {%
duinkerken2006comparison}%
\begin{APACrefauthors}%
Duinkerken, M\BPBI B.%
, Ottjes, J\BPBI A.%
\BCBL {}\ \BBA {} Lodewijks, G.%
\end{APACrefauthors}%
\unskip\
\newblock
\APACrefYearMonthDay{2006}{}{}.
\newblock
{\BBOQ}\APACrefatitle {Comparison of Routing Strategies for AGV Systems using
  Simulation} {Comparison of routing strategies for agv systems using
  simulation}.{\BBCQ}
\newblock
\BIn{} \APACrefbtitle {Proceedings of the 2006 Winter Simulation Conference}
  {Proceedings of the 2006 winter simulation conference}\ (\BPGS\ 1523--1530).
\PrintBackRefs{\CurrentBib}

\bibitem [\protect \citeauthoryear {%
Guillaume%
, Michael%
, Abdelghani%
, Damien%
\BCBL {}\ \BBA {} Mohamed%
}{%
Guillaume%
\ \protect \BOthers {.}}{%
{\protect \APACyear {2017}}%
}]{%
guillaume2017}
\APACinsertmetastar {%
guillaume2017}%
\begin{APACrefauthors}%
Guillaume, D.%
, Michael, D.%
, Abdelghani, B.%
, Damien, T.%
\BCBL {}\ \BBA {} Mohamed, D.%
\end{APACrefauthors}%
\unskip\
\newblock
\APACrefYearMonthDay{2017}{}{}.
\newblock
{\BBOQ}\APACrefatitle {Decentralized Motion Planning and Scheduling of AGVs in
  FMS} {Decentralized motion planning and scheduling of agvs in fms}.{\BBCQ}
\newblock
\APACjournalVolNumPages{Transactions on Industrial Informatics}{}{}{}.
\PrintBackRefs{\CurrentBib}

\bibitem [\protect \citeauthoryear {%
Gurobi~Optimization%
}{%
Gurobi~Optimization%
}{%
{\protect \APACyear {2021}}%
}]{%
gurobi}
\APACinsertmetastar {%
gurobi}%
\begin{APACrefauthors}%
Gurobi~Optimization, L.%
\end{APACrefauthors}%
\unskip\
\newblock
\APACrefYearMonthDay{2021}{}{}.
\newblock
\APACrefbtitle {Gurobi Optimizer Reference Manual.} {Gurobi optimizer reference
  manual.}
\newblock
\begin{APACrefURL} \url{http://www.gurobi.com} \end{APACrefURL}
\PrintBackRefs{\CurrentBib}

\bibitem [\protect \citeauthoryear {%
Herrero-Perez%
\ \BBA {} Martinez-Barbera%
}{%
Herrero-Perez%
\ \BBA {} Martinez-Barbera%
}{%
{\protect \APACyear {2010}}%
}]{%
herrero2010}
\APACinsertmetastar {%
herrero2010}%
\begin{APACrefauthors}%
Herrero-Perez, D.%
\BCBT {}\ \BBA {} Martinez-Barbera, H.%
\end{APACrefauthors}%
\unskip\
\newblock
\APACrefYearMonthDay{2010}{}{}.
\newblock
{\BBOQ}\APACrefatitle {Modeling Distributed Transportation Systems Composed of
  Flexible Automated Guided Vehicles in Flexible Manufacturing Systems}
  {Modeling distributed transportation systems composed of flexible automated
  guided vehicles in flexible manufacturing systems}.{\BBCQ}
\newblock
\APACjournalVolNumPages{IEEE Transactions on Industrial
  Informatics}{6}{2}{166--180}.
\PrintBackRefs{\CurrentBib}

\bibitem [\protect \citeauthoryear {%
Korsah%
, Stentz%
\BCBL {}\ \BBA {} Dias%
}{%
Korsah%
\ \protect \BOthers {.}}{%
{\protect \APACyear {2013}}%
}]{%
korsah2013}
\APACinsertmetastar {%
korsah2013}%
\begin{APACrefauthors}%
Korsah, G\BPBI A.%
, Stentz, A.%
\BCBL {}\ \BBA {} Dias, M\BPBI B.%
\end{APACrefauthors}%
\unskip\
\newblock
\APACrefYearMonthDay{2013}{}{}.
\newblock
{\BBOQ}\APACrefatitle {A Comprehensive Taxonomy for Multi-Robot Task
  Allocation} {A comprehensive taxonomy for multi-robot task
  allocation}.{\BBCQ}
\newblock
\APACjournalVolNumPages{The International Journal of Robotics
  Research}{32}{12}{1495--1512}.
\PrintBackRefs{\CurrentBib}

\bibitem [\protect \citeauthoryear {%
Lee%
, Zaheer%
\BCBL {}\ \BBA {} Kim%
}{%
Lee%
\ \protect \BOthers {.}}{%
{\protect \APACyear {2014}}%
}]{%
lee2014}
\APACinsertmetastar {%
lee2014}%
\begin{APACrefauthors}%
Lee, D\BHBI H.%
, Zaheer, S\BPBI A.%
\BCBL {}\ \BBA {} Kim, J\BHBI H.%
\end{APACrefauthors}%
\unskip\
\newblock
\APACrefYearMonthDay{2014}{}{}.
\newblock
{\BBOQ}\APACrefatitle {A Resource-oriented, Decentralized Auction algorithm for
  Multirobot Task Allocation} {A resource-oriented, decentralized auction
  algorithm for multirobot task allocation}.{\BBCQ}
\newblock
\APACjournalVolNumPages{IEEE Transactions on Automation Science and
  Engineering}{12}{4}{1469--1481}.
\PrintBackRefs{\CurrentBib}

\bibitem [\protect \citeauthoryear {%
Miller%
, Tucker%
\BCBL {}\ \BBA {} Zemlin%
}{%
Miller%
\ \protect \BOthers {.}}{%
{\protect \APACyear {1960}}%
}]{%
miller1960}
\APACinsertmetastar {%
miller1960}%
\begin{APACrefauthors}%
Miller, C\BPBI E.%
, Tucker, A\BPBI W.%
\BCBL {}\ \BBA {} Zemlin, R\BPBI A.%
\end{APACrefauthors}%
\unskip\
\newblock
\APACrefYearMonthDay{1960}{}{}.
\newblock
{\BBOQ}\APACrefatitle {Integer Programming Formulation of Traveling Salesman
  Problems} {Integer programming formulation of traveling salesman
  problems}.{\BBCQ}
\newblock
\APACjournalVolNumPages{Journal of the ACM}{7}{4}{326--329}.
\PrintBackRefs{\CurrentBib}

\bibitem [\protect \citeauthoryear {%
Miyamoto%
\ \BBA {} Inoue%
}{%
Miyamoto%
\ \BBA {} Inoue%
}{%
{\protect \APACyear {2016}}%
}]{%
miyamoto2016local}
\APACinsertmetastar {%
miyamoto2016local}%
\begin{APACrefauthors}%
Miyamoto, T.%
\BCBT {}\ \BBA {} Inoue, K.%
\end{APACrefauthors}%
\unskip\
\newblock
\APACrefYearMonthDay{2016}{}{}.
\newblock
{\BBOQ}\APACrefatitle {Local and Random Searches for Dispatch and Conflict-free
  Routing Problem of Capacitated AGV Systems} {Local and random searches for
  dispatch and conflict-free routing problem of capacitated agv
  systems}.{\BBCQ}
\newblock
\APACjournalVolNumPages{Computers \& Industrial Engineering}{91}{}{1--9}.
\PrintBackRefs{\CurrentBib}

\bibitem [\protect \citeauthoryear {%
Reinelt%
}{%
Reinelt%
}{%
{\protect \APACyear {1991}}%
}]{%
reinelt1991}
\APACinsertmetastar {%
reinelt1991}%
\begin{APACrefauthors}%
Reinelt, G.%
\end{APACrefauthors}%
\unskip\
\newblock
\APACrefYearMonthDay{1991}{}{}.
\newblock
{\BBOQ}\APACrefatitle {TSPLIB—A Traveling Salesman Problem Library}
  {Tsplib—a traveling salesman problem library}.{\BBCQ}
\newblock
\APACjournalVolNumPages{ORSA Journal on Computing}{3}{4}{376--384}.
\PrintBackRefs{\CurrentBib}

\bibitem [\protect \citeauthoryear {%
Shi%
\ \BBA {} Ng%
}{%
Shi%
\ \BBA {} Ng%
}{%
{\protect \APACyear {2018}}%
}]{%
shi2018collision}
\APACinsertmetastar {%
shi2018collision}%
\begin{APACrefauthors}%
Shi, Z.%
\BCBT {}\ \BBA {} Ng, W\BPBI K.%
\end{APACrefauthors}%
\unskip\
\newblock
\APACrefYearMonthDay{2018}{}{}.
\newblock
{\BBOQ}\APACrefatitle {A Collision-free Path Planning Algorithm for Unmanned
  Aerial Vehicle Delivery} {A collision-free path planning algorithm for
  unmanned aerial vehicle delivery}.{\BBCQ}
\newblock
\BIn{} \APACrefbtitle {2018 International Conference on Unmanned Aircraft
  Systems (ICUAS)} {2018 international conference on unmanned aircraft systems
  (icuas)}\ (\BPGS\ 358--362).
\PrintBackRefs{\CurrentBib}

\bibitem [\protect \citeauthoryear {%
Smolic-Rocak%
, Bogdan%
, Kovacic%
\BCBL {}\ \BBA {} Petrovic%
}{%
Smolic-Rocak%
\ \protect \BOthers {.}}{%
{\protect \APACyear {2009}}%
}]{%
smolic2009time}
\APACinsertmetastar {%
smolic2009time}%
\begin{APACrefauthors}%
Smolic-Rocak, N.%
, Bogdan, S.%
, Kovacic, Z.%
\BCBL {}\ \BBA {} Petrovic, T.%
\end{APACrefauthors}%
\unskip\
\newblock
\APACrefYearMonthDay{2009}{}{}.
\newblock
{\BBOQ}\APACrefatitle {Time Windows based Dynamic Routing in Multi-AGV Systems}
  {Time windows based dynamic routing in multi-agv systems}.{\BBCQ}
\newblock
\APACjournalVolNumPages{IEEE Transactions on Automation Science and
  Engineering}{7}{1}{151--155}.
\PrintBackRefs{\CurrentBib}

\bibitem [\protect \citeauthoryear {%
Spensieri%
, Carlson%
, Ekstedt%
\BCBL {}\ \BBA {} Bohlin%
}{%
Spensieri%
\ \protect \BOthers {.}}{%
{\protect \APACyear {2015}}%
}]{%
spensieri2015}
\APACinsertmetastar {%
spensieri2015}%
\begin{APACrefauthors}%
Spensieri, D.%
, Carlson, J\BPBI S.%
, Ekstedt, F.%
\BCBL {}\ \BBA {} Bohlin, R.%
\end{APACrefauthors}%
\unskip\
\newblock
\APACrefYearMonthDay{2015}{}{}.
\newblock
{\BBOQ}\APACrefatitle {An Iterative Approach for Collision Free Routing and
  Scheduling in Multirobot Stations} {An iterative approach for collision free
  routing and scheduling in multirobot stations}.{\BBCQ}
\newblock
\APACjournalVolNumPages{IEEE Transactions on Automation Science and
  Engineering}{13}{2}{950--962}.
\PrintBackRefs{\CurrentBib}

\bibitem [\protect \citeauthoryear {%
Toth%
\ \BBA {} Vigo%
}{%
Toth%
\ \BBA {} Vigo%
}{%
{\protect \APACyear {2014}}%
}]{%
toth2014}
\APACinsertmetastar {%
toth2014}%
\begin{APACrefauthors}%
Toth, P.%
\BCBT {}\ \BBA {} Vigo, D.%
\end{APACrefauthors}%
\unskip\
\newblock
\APACrefYear{2014}.
\newblock
\APACrefbtitle {Vehicle Routing: Problems, Methods, and Applications} {Vehicle
  routing: Problems, methods, and applications}.
\newblock
\APACaddressPublisher{}{SIAM}.
\PrintBackRefs{\CurrentBib}

\bibitem [\protect \citeauthoryear {%
Xin%
\ \protect \BOthers {.}}{%
Xin%
\ \protect \BOthers {.}}{%
{\protect \APACyear {2020}}%
}]{%
xin2020time}
\APACinsertmetastar {%
xin2020time}%
\begin{APACrefauthors}%
Xin, J.%
, Meng, C.%
, Schulte, F.%
, Peng, J.%
, Liu, Y.%
\BCBL {}\ \BBA {} Negenborn, R\BPBI R.%
\end{APACrefauthors}%
\unskip\
\newblock
\APACrefYearMonthDay{2020}{}{}.
\newblock
{\BBOQ}\APACrefatitle {A Time-Space Network Model for Collision-Free Routing of
  Planar Motions in a Multirobot Station} {A time-space network model for
  collision-free routing of planar motions in a multirobot station}.{\BBCQ}
\newblock
\APACjournalVolNumPages{IEEE Transactions on Industrial
  Informatics}{16}{10}{6413--6422}.
\PrintBackRefs{\CurrentBib}

\bibitem [\protect \citeauthoryear {%
Xin%
, Negenborn%
, Corman%
\BCBL {}\ \BBA {} Lodewijks%
}{%
Xin%
\ \protect \BOthers {.}}{%
{\protect \APACyear {2015}}%
}]{%
xin2015control}
\APACinsertmetastar {%
xin2015control}%
\begin{APACrefauthors}%
Xin, J.%
, Negenborn, R\BPBI R.%
, Corman, F.%
\BCBL {}\ \BBA {} Lodewijks, G.%
\end{APACrefauthors}%
\unskip\
\newblock
\APACrefYearMonthDay{2015}{}{}.
\newblock
{\BBOQ}\APACrefatitle {Control of Interacting Machines in Automated Container
  Terminals using a Sequential Planning Approach for Collision Avoidance}
  {Control of interacting machines in automated container terminals using a
  sequential planning approach for collision avoidance}.{\BBCQ}
\newblock
\APACjournalVolNumPages{Transportation Research Part C: Emerging
  Technologies}{60}{}{377--396}.
\PrintBackRefs{\CurrentBib}

\end{thebibliography}

\end{document}